\documentclass[11pt]{article}
\usepackage{pifont}

\usepackage{amsfonts}
\usepackage{latexsym}
\usepackage{amsmath}
\usepackage{amssymb}
\usepackage{color}
\newtheorem{theorem}{\bf{Theorem}}[section]
\newtheorem{definition}[theorem]{\bf{Definition}}
\newtheorem{example}[theorem]{\bf{Example}}
\newtheorem{corol}[theorem]{\bf{Corollary}}
\newtheorem{proposition}[theorem]{\bf{Proposition}}
\newtheorem{remark}[theorem]{\bf{Remark}}

\begin{document}
\hoffset = -2truecm \voffset = -2truecm
\title{\hspace{-0.9cm}\bf A new comparison theorem of multidimensional\footnotetext{{\small Supported by the National Basic Research Program of China (No. 2007CB814901), and WCU (World Class
University) program of the Korea Science and Engineering Foundation (R31-20007).}
} \\
\hspace{-13.2cm}BSDEs}
\date{}

\author{\hspace{-8.7cm} \bf{Pan-yu Wu$^{1,*}$, Zeng-jing Chen$^{1,2}$}\\
\hspace{-4.9cm}
 $^{1}${\small School of Mathematics, Shandong
University, Jinan, 250100, China }\\
\hspace{-9.2cm}
$^{2}${\small Ajou University, Suwon, 443749, Korea }\\
\hspace{-9.8cm}
\small{(E-mail:
$^*$panyuaza@yahoo.com.cn)}}
\maketitle \vspace{.1in}

\hspace{-0.7cm} {\small\bf Abstract\,\,\,\,} {\small In this paper,
we first study
a property about the generator $g$ of Backward Stochastic Differential
Equation (BSDE) when the price of contingent
claims can be represented by a multidimensional BSDE in the no-arbitrage financial market.
Furthermore, motivated by the
behavior of agent in finance market, we introduce a new total order $\succsim^q$ on $\mathbb{R}^n$ and obtain a
necessary and sufficient
condition for comparison theorem of multidimensional BSDEs under this order. We also give some further
results for special order $\succsim^q$.}
\\
\\
{\small\bf Keywords\,\,\,\,} {\small Backward stochastic differential equation, Comparison theorem, Viability property}\\
{\small\bf 2000 MR Subject Classification}\hspace{0.5cm} {\small 60H10; 60H30}\\
\section{ Preliminaries}
In this section, we shall introduce some notations and assumptions
which are needed in the following analysis.

Let $(\Omega,\mathcal{F},P)$ be a probability space and $(W_t)_{t\geq0}$ a
standard $d$-dimensional Brownian motion defined on this probability
space. Furthermore, let $(\mathcal{F}_t)_{t\geq0}$ be the filtration
generated by the Brownian motion $(W_t)_{t\geq0}$, that is
$\mathcal{F}_t=\sigma(W_s;s\leq t)$. We define the usual $P$-augmentation to
each $\mathcal{F}_t$ such that $(\mathcal{F}_t)_{t\geq0}$ is right continuous and
complete. We denote by $\mathbb{R}^n$ the $n$-dimensional Euclidean space,
equipped with the standard inner product $\langle\cdot,\cdot\rangle$ and
the Euclidean norm $|\cdot|$. We also denote by $\mathbb{R}^{n\times d}$ the
collection of all $n\times d$ real matrices, and for matrix
$z=(z_{ij})_{n\times d}$, we denote $z_i:=(z_{i1},\cdots,z_{id})^{\mathbf{T}}$
and $|z|:=\sqrt{tr(zz^{\mathbf{T}})}$, where $z^{\mathbf{T}}$ represents the transpose of $z$. Now, we define the following usual
spaces of random variables or processes:
\\$\bullet\ \ L^2(\Omega,\mathcal{F}_t,P;K)=\{\xi|\ \xi$ is $K$-valued $\mathcal{F}_t$-measurable random variable and $\mathbb{E}[\xi^2]<\infty\}$, where $K$ is a subset of $\mathbb{R}^n$;
\\$\bullet\ \ \mathcal{S}^2_T=\{\psi|\ \psi$ is $\mathbb{R}^n$-valued progressively measurable process and $\mathbb{E}[\sup_{0\leq t\leq T}|\psi_t|^2]<\infty\}$;
\\$\bullet\ \ \mathcal{H}^2_T=\{\psi|\ \psi$ is $\mathbb{R}^{n\times d}$-valued progressively measurable process and $\mathbb{E}[\int_0^T|\psi_t|^2dt]<\infty\}.$
\par
Consider the following Backward Stochastic Differential Equation (BSDE for short):
\begin{equation}\label{eq11}
Y_t=\xi+\int_t^T g(s,Y_s,Z_s)ds-\int_t^T Z_sdW_s,\quad 0\leq t\leq T,
\end{equation}
where $\xi$ is a given $n$-dimensional random variable, and $g$ is called the generator of BSDE (\ref{eq11}), which is defined as $\Omega\times[0,T]\times\mathbb{R}^n\times\mathbb{R}^{n\times d}\to\mathbb{R}^n$, such that the process $(g(t,y,z))_{t\in[0,T]}$ is progressively measurable for each $(y,z)$ in $\mathbb{R}^n\times\mathbb{R}^{n\times d}$. We make the following assumptions (A1)-(A4) throughout this paper:
\\$\bullet\ \ $(A1) For any $(y,z)\in\mathbb{R}^n\times\mathbb{R}^{n\times d}$, $t\rightarrow g(t,y,z)$ is continuous, $P$-a.s.;
\\$\bullet\ \ $(A2) There exists a constant $\mu>0,$ such that for any $t\in[0,T]$, $(y,z),(y',z')\in\mathbb{R}^n\times\mathbb{R}^{n\times d},$ we have
\begin{center}$|g(t,y,z)-g(t,y',z')|\leq \mu(|y-y'|+|z-z'|),$ \ \ $P$-a.s.;\end{center}
$\bullet\ \ $(A3) $g(t,0,0)\in \mathcal{S}^2_T;$
\\$\bullet\ \ $(A4) $\xi\in L^2(\Omega,\mathcal{F}_T,P;\mathbb{R}^n).$
\par
Under the assumptions (A1)-(A4), Pardoux and Peng \cite{par1} showed
that  BSDE (\ref{eq11}) has a unique adapted solution $(Y,Z)$
belonging to $\mathcal{S}^2_T\times\mathcal{H}^2_T$. The viability property of a
stochastic process as a classical notion in stochastic context was
first discussed in \cite{a}.  Buckdahn et al. \cite{bk} studied the
viability property of BSDE (\ref{eq11}). In Theorem 2.5 of \cite{bk},
they obtained the Backward Stochastic Viability Property (BSVP for
short) for BSDE. In next section, using the BSVP, we
study the property about the generator $g$ of BSDE when the price
vector of contingent claims can be represented by a multidimensional
BSDE in the no-arbitrage financial market. Some other applications of BSDE
in financial mathematics can be found for example in [3, 5, 6,
8].
\par
The comparison theorem and related converse comparison theorem for one
dimensional BSDEs were important results in the theory of BSDE
first due to Peng \cite{p92} and Coquet et al. \cite{f} respectively, and later
generalized by several authors (for example see \cite{el}). Combining
comparison theorem and converse comparison theorem for
one-multidimensional BSDEs, we can get the following result:
For any $0\leq u\leq T$, consider the following two BSDEs,
\begin{equation}\label{eq0}
Y_t^i=\xi^i+\int_t^{u} g^i(s,Y_s^i,Z_s^i)ds-\int_t^{u} Z_s^idW_s,\quad 0\leq t\leq u,\quad i=1,2,
\end{equation}
where $g^1$ and $g^2$ satisfy (A1)-(A3), then the statements (i$'$) and (ii$'$) are equivalent:\\
$\bullet\ $ (i$'$) For any $u\in[0,T]$, $\xi^1,\xi^2\in L^2(\Omega,\mathcal{F}_{u},P;\mathbb{R})$, if $\xi^1\geq\xi^2$, then $Y_t^1\geq Y_t^2$, $P$-a.s., for $t\in[0,u]$.
\\$\bullet\ $ (ii$'$) For any $ t\in[0,T]$, $(y,z)\in\mathbb{R}\times\mathbb{R}^d$, $g^1(t,y,z)\geq g^2(t,y,z)$, $P$-a.s..
\par
A natural question is whether the above  equivalence still holds for
multidimensional BSDEs. To answer this question, the key point is
how to define ``order'' or ``preference'' on $\mathbb{R}^n$. Hu
and Peng \cite{hu} considered the case where  the ``order" $y^1\geq
y^2$ on $\mathbb{R}^n$ is in the sense of $y^1_i\geq y^2_i,$ for all
$i=1,2,\cdots n$, where $y_i^1$ and $y_i^2$
are the $i$-th component of $y^1$ and $y^2$ respectively. They obtained a necessary and sufficient condition
of (i$'$) for multidimensional BSDEs. However,
in finance, such a preference is not enough to
describe the behavior of agents. For example, in financial market, let $y^1=(y_1^1, \cdots, y_n^1)$, $y^2=(y_1^2, \cdots, y_n^2)$ be
two portfolios consisting of $n$-basic contingent claims, $q$ be the price vector of those contingent claims.
Agents often like to compare the value $\langle y^1,q\rangle$ and $\langle y^2,q\rangle$ of
portfolios. In this case, it is natural to define a total order
$\succsim^q$ on $\mathbb{R}^n$ via $q$. What is comparison theorem under this order? In Section 3, we re-state (i$'$) for
multidimensional BSDEs, and obtain a necessary and sufficient
condition for comparison theorem of multidimensional BSDEs under the
new total order. The result is another application of BSVP. We also give some further results for special total order $\succsim^q$.

\section{BSVP and its Application}
Let us recall the definition of BSVP from \cite{bk}.
\begin{definition}\label{de13}
Let $K$ be a nonempty, convex closed set in $\mathbb{R}^n$. Then we call the BSDE (\ref{eq11}) enjoys the BSVP in $K$ if: for any $u\in[0,T ]$, $\xi\in L^2(\Omega,\mathcal{F}_{u},P;K)$, the unique solution $(Y,Z)\in\mathcal{S}^2_{u}\times\mathcal{H}^2_{u}$ to the BSDE (\ref{eq11}) over time interval $[0,u]$, given by
\begin{equation}\label{eq12}
Y_t=\xi+\int_t^{u} g(s,Y_s,Z_s)ds-\int_t^{u} Z_sdW_s,\quad 0\leq t\leq u,
\end{equation}
satisfies for any $t\in[0,u]$, $Y_t\in K$, $P$-a.s..
\end{definition}
For completeness, we recall the necessary and sufficient condition
of BSVP for the BSDEs in \cite{bk}.
\begin{proposition}\label{th14}
Let $K$ be a nonempty, convex closed set in $\mathbb{R}^n$. Suppose that $g$
satisfies (A1)-(A3). Then  BSDE (\ref{eq11}) enjoys the BSVP in $K$
if and only if for any $(t,z)\in[0,T]\times\mathbb{R}^{n\times d}$ and any $ y\in\mathbb{R}^n$
such that $d_K^2(\cdot)$ is twice differentiable at $y$,
\begin{equation}\label{eq13}
4\left\langle y-\Pi_K(y),g(t,\Pi_K(y),z)\right\rangle\leq\langle D^2d^2_K(y)z,z\rangle+Cd_K^2(y),\quad P\hbox{-a.s.},
\end{equation}
where $C >0$ is a constant which does not depend on $(t, y, z)$, $\Pi_K(y)$ is the projection of $y$ onto $K$, $d_K(y)$ represents the distance between $y$ and $K$.
\end{proposition}
\par
Now we give an application of BSVP in the no-arbitrage financial
market. We get a property of $g$ when the price vector of contingent
claims can be represented by a multidimensional BSDE (see \cite{el}
for details). Suppose that there are $n$ kinds of contingent claims
in the market. The price vector of this $n$ contingent claims is a random process $(Y_t)_{0\leq t\leq T}$ with
$Y_t\in L^2(\Omega,\mathcal{F}_{t},P;\mathbb{R}^n)$. Let $X^1,X^2$ be two different kinds of risk
positions. Furthermore, without loss of generality, assume that the
terminal value of $X^1$ is bigger than that of $X^2$, in other words,
$\langle Y_T,q\rangle\geq0$, where $q:=(X^1-X^2)/|X^1-X^2|$. Then we
have the following theorem.
\begin{theorem}\label{pr15}
If
\begin{equation}\label{eq14}
Y_t=Y_T+\int_t^T g(s,Y_s,Z_s)ds-\int_t^T Z_sdW_s,\quad 0\leq t\leq T,
\end{equation}
then
$$-4\langle y,q\rangle^-\left\langle q,g(t,y+\langle y,q\rangle^-q,z)\right\rangle\leq2I_{\langle y,q\rangle<0}\sum_{i,j=1}^nq_iq_jz_i^{\mathbf{T}}z_j+C(\langle y,q\rangle^-)^2,  \quad P\hbox{-a.s.},$$where $C >0$ is a constant which does not depend on $(t, y, z)$.
\end{theorem}
{\it Proof.}\,\, Since the market is no-arbitrage, we have if
$\langle Y_T,q\rangle\geq0$, then $\langle Y_t,q\rangle\geq0,$
$P$-a.s, for any $t\in[0,T]$ (for details see \cite{el}). That is,
BSDE (\ref{eq14}) enjoys the BSVP in $K$, where
$K:=\{y\in\mathbb{R}^n|\ \langle y,q\rangle\geq 0\}$. Clearly, $K$ is a
nonempty, convex closed set of $\mathbb{R}^n$.
\par
If $y\in K$, then $\Pi_K(y)=y,d_K(y)=0.$ If $y\not\in K$, we have
\begin{displaymath}
\left\{ \begin{array}{ll}
\langle \Pi_K(y),q\rangle=0 \\
\Pi_K(y)-d_K(y)q=y
\end{array} \right..
\end{displaymath}
Let $\Pi_K(y)=(u_1,u_2,\dots,u_n)$, solving the above two equations, we get
\begin{displaymath}
\left( \begin{array}{c}
u_1\\
u_2\\
\vdots\\
u_n\\
d_K(y)
\end{array} \right)=
\left( \begin{array}{ccccc}
1-q_1^2 & -q_1q_2 & \cdots & -q_1q_n & q_1 \\
-q_1q_2 & 1-q_2^2 & \cdots & -q_2q_n & q_2\\
\cdots&\cdots&\cdots&\cdots&\cdots\\
-q_1q_n & -q_2q_n &\cdots&1-q_n^2&q_n\\
-q_1&-q_2&\cdots&-q_n&1
\end{array} \right)
\left( \begin{array}{c}
y_1\\
y_2\\
\vdots\\
y_n\\
0
\end{array} \right),
\end{displaymath}
that is, $\Pi_K(y)=y-\langle y,q\rangle q,$ $d_K(y)=-\langle y,q\rangle$. Consequently, for any $y\in\mathbb{R}^n$, we get $\Pi_K(y)=y+\langle y,q\rangle^- q$, $d_K(y)=\langle y,q\rangle^-$. Therefore, for $y\in\mathbb{R}^n$, $D^2d^2_K(y)=2qq^{\mathbf{T}}I_{\langle y,q\rangle<0}$.
Due to Proposition \ref{th14}, $g$ must satisfies $$-4\langle y,q\rangle^-\left\langle q,g(t,y+\langle y,q\rangle^-q,z)\right\rangle\leq2I_{\langle y,q\rangle<0}\sum_{i,j=1}^nq_iq_jz_i^{\mathbf{T}}z_j+
C(\langle y,q\rangle^-)^2,  \quad P\hbox{-a.s.}.$$
The proof of Theorem \ref{pr15} is completed. $\hfill{} \Box$
\section{Comparison Theorem for Multidimensional BSDEs under the Total Order $\succsim^q$}
In this section, we first introduce the definition of a total order
on $\mathbb{R}^n$ denoted by $\succsim^q$  and then prove the comparison
theorem for multidimensional BSDEs under this total order.
\begin{definition}\label{de11}
Let $q\in\mathbb{R}^n$ be any fixed nonvanishing vector. For any $y^1,y^2\in \mathbb{R}^n,$ we call $y^1$ bigger (or better) than $y^2$ under $q$, denote $y^1\succsim^q y^2$, if $\langle y^1,q\rangle\geq \langle y^2,q\rangle$.
\end{definition}
\begin{remark}\label{re12}
$(1)$ Obviously, $y^1\succsim^q y^2$ if and only if $y^1\succsim^{q/|q|} y^2$. So without loss of generality, we assume $q$ be a unit vector in the sequel.
\\$(2)$ $\succsim^q$ is a total order on $\mathbb{R}^n$, which can be used to compare any two elements in $\mathbb{R}^n$.
\end{remark}
\par
We now begin to prove the comparison theorem of multidimensional
BSDEs under the total order $\succsim^q$.
\par
For any $0\leq u\leq T$, consider the following two BSDEs,
\begin{equation}\label{eq21}
Y_t^i=\xi^i+\int_t^{u} g^i(s,Y_s^i,Z_s^i)ds-\int_t^{u} Z_s^idW_s,\quad 0\leq t\leq u,\quad i=1,2.
\end{equation}
\begin{theorem}\label{th21}
Suppose that $g^1$ and $g^2$ satisfy (A1)-(A3). Then, the following two statements (i) and (ii) are equivalent.
\\$\bullet$ (i) For any $u\in[0,T]$, $\xi^1,\xi^2\in L^2(\Omega,\mathcal{F}_{u},P;\mathbb{R}^n)$ such that $\xi^1\succsim^q\xi^2 $, then the unique solutions $(Y^1,Z^1)$ and $(Y^2,Z^2)$ in $\mathcal{S}^2_{u}\times\mathcal{H}^2_{u}$ to BSDEs (\ref{eq21}) over time interval $[0,u]$ satisfy $$Y_t^1\succsim^q Y_t^2, \quad P\hbox{-a.s.},\quad \forall t\in[0,u].$$
$\bullet$ (ii) For any $ t\in[0,T]$, $(y,z),(y',z')\in\mathbb{R}^n\times\mathbb{R}^{n\times d}$, we have
\begin{equation*}\label{eq22}
-4\langle y,q\rangle^-\left\langle q,\,g^1(t,y+\langle y,q\rangle^-q+y',z)-g^2(t,y',z')\right\rangle
\end{equation*}
\begin{equation}\label{eq22}
\leq2I_{\langle y,q\rangle<0}\sum_{i,j=1}^nq_iq_j(z_i-z'_i)^{\mathbf{T}}(z_j-z'_j)+C(\langle y,q\rangle^-)^2, \quad P\hbox{-a.s.},
\end{equation}
where $C >0$ is a constant independent of $(t, y, z)$.
\end{theorem}
{\it Proof.}\,\, It is obviously that (\ref{eq21})$\Leftrightarrow$
\begin{displaymath}
\left\{ \begin{array}{c}
Y_t^1-Y_t^2=\xi^1-\xi^2+\int_t^{u} [g^1(s,Y_s^1,Z_s^1)-g^2(s,Y^2_s,Z_s^2)]ds-\int_t^{u} (Z_s^1-Z_s^2)dW_s,\ \ \ 0\leq t\leq u,\\
Y_t^2=\xi^2+\int_t^{u} g^2(s,Y_s^2,Z_s^2)ds-\int_t^{u} Z_s^2dW_s,\ \ \ 0\leq t\leq u.\qquad\qquad\qquad\qquad\qquad\qquad\qquad\ \
\end{array} \right.
\end{displaymath}
\par
Let\begin{displaymath}
\bar{Y}_t=
\left(\begin{array}{c}
\bar{Y}_t^1\\ \bar{Y}_t^2
\end{array} \right):=\left(\begin{array}{c}
Y_t^1-Y_t^2\\ Y_t^2
\end{array} \right),
\end{displaymath}
\begin{displaymath}
\bar{Z}_t=
\left(\begin{array}{c}
\bar{Z}_t^1\\ \bar{Z}_t^2
\end{array} \right):=\left(\begin{array}{c}
Z_t^1-Z_t^2\\ Z_t^2
\end{array} \right),\
\bar{\xi}=
\left(\begin{array}{c}
\bar{\xi}^1\\ \bar{\xi}^2
\end{array} \right):=\left(\begin{array}{c}
\xi^1-\xi^2\\\xi^2
\end{array} \right).
\end{displaymath}
Then (i) is equivalent to the following statement (iii):\\
$\bullet$ (iii) For any $ u\in[0,T ]$, $\bar{\xi}=
\left(\begin{array}{c}
\bar{\xi}^1\\ \bar{\xi}^2
\end{array} \right)\in L^2(\Omega,\mathcal{F}_{u},P;\mathbb{R}^{2n}) $
such that $\bar{\xi}^1\succsim^q0$, the unique solution $(\bar{Y},\bar{Z})$ to the following $2n$-dimensional BSDE (\ref{eq23}) over time interval $[0,u]$ satisfies $\bar{Y}^1_t\succsim^q0,$ $P$-a.s.,
\begin{equation}\label{eq23}
\bar{Y}_t=\bar{\xi}+\int_t^{u} \bar{g}(s,\bar{Y}_s,\bar{Z}_s)ds-\int_t^{u} \bar{Z}_sdW_s,\quad 0\leq t\leq u,
\end{equation}
where for
$y=\left(\begin{array}{c}y^1\\y^2\end{array} \right)$,
$z=\left(\begin{array}{c}z^1\\ z^2\end{array} \right)$, we have
$$\bar{g}(t,y,z)=\left(\begin{array}{c}g^1(t,y^1+y^2,z^1+z^2)-g^2(t,y^2,z^2)\\ g^2(t,y^2,z^2)\end{array} \right).$$
\par
The statement (iii) means that the BSDE (\ref{eq23}) satisfies BSVP
in $K:=\{x\in\mathbb{R}^n|\ x\succsim^q0\}\times\mathbb{R}^n$. Obviously, $K$ is a
nonempty, convex closed subset of $\mathbb{R}^{2n}$. Similarly to the proof of Theorem
\ref{pr15}, we have
$$\Pi_K(y)=\left(\begin{array}{c}y^1+\langle y^1,q\rangle^- q\\y^2\end{array} \right),\  d_K^2(y)=(\langle y^1,q\rangle^-)^2,\ D^2d^2_K(y)=\left(\begin{array}{cc}2qq^{\mathbf{T}}I_{\langle y^1,q\rangle<0}&\mathbf{0}\\ \mathbf{0}&\mathbf{0}\end{array} \right),$$
 where $\mathbf{0}$ is the $n$-order zero matrix. Applying Proposition \ref{th14}, we obtain that the statement (iii) is equivalent to:
\begin{equation*}\textstyle -4\langle y^1,q\rangle^-\left\langle q,\,g^1(t,y^1+\langle y^1,q\rangle^- q+y^2,z^1+z^2)-g^2(t,y^2,z^2)\right\rangle
\end{equation*}
\begin{equation}\label{3}
\leq2I_{\langle y^1,q\rangle<0}\sum_{i,j=1}^nq_iq_j(z_i^1)^{\mathbf{T}}z_j^1+C(\langle y^1,q\rangle^-)^2,\quad P\hbox{-a.s.},
\end{equation}
where $C >0$ is a constant independent of $(t, y, z)$.
Let$$\left(\begin{array}{c}y^1\\y^2\end{array} \right)=\left(\begin{array}{c}y\\y'\end{array} \right),
\left(\begin{array}{c}z^1\\z^2\end{array} \right)=\left(\begin{array}{c}z-z'\\z'\end{array} \right),$$
it is then clear that the above inequality (\ref{3}) becomes (\ref{eq22}). The proof of Theorem \ref{th21} is completed. $\hfill{} \Box$
\par
Furthermore, we can get the following Theorem \ref{th22}.
\begin{theorem}\label{th22}
If (i) holds, then for any $(t,y,z)\in[0,T]\times\mathbb{R}^n\times\mathbb{R}^{n\times d}$, we have $g^1(t,y,z)\succsim^qg^2(t,y,z)$, $P$-a.s..
\end{theorem}
{\it Proof.}\,\, Because $q=(q_1,q_2,\cdots,q_n)^{\mathbf{T}}\neq0$, there exists $ q_i\neq0$ for some $i$. Without loss of generality, we suppose $q_i>0$. Set $y^k=-{1\over k}e^i$, where the components of $e^i$ are 0 except the $i$th component which is $1$, and $k$ is an arbitrary number in $\mathbb{N}^*$. Then, $$\langle y^k,q\rangle^-={1\over k}q_i,\ \ \ \ y^k+\langle y^k,q\rangle^- q={1\over k}[q_iq-e^i].$$
\par
So for any $(t,y,z)\in[0,T]\times\mathbb{R}^n\times\mathbb{R}^{n\times d}$, substituting $(y^k,z)$ and $(y,z)$ into $(y,z)$ and $(y',z')$ in inequality (\ref{eq22}) respectively, we have
$$\left\langle q,\,g^1(t,-\frac 2 k e^i+\frac 1 k q_iq+y,z)-g^2(t,y,z)\right\rangle\geq-\frac{Cq_i}{4k},\quad P\hbox{-a.s.},$$
where $C>0$ is a constant which does not depend on $(t,y,z).$
As $k\to\infty$, it follows from $(A2)$ that $g^1(t,y,z)\succsim^qg^2(t,y,z)$, $P$-a.s.. The proof of Theorem \ref{th22} is completed.$\hfill{} \Box$
\par
Note that the converse of Theorem \ref{th22} is not true (see Example \ref{ex1} given latter). As an application of Theorem \ref{th22} we immediately have the following corollary.
\begin{corol}\label{co23}
If (i) holds for both $q=e^i$, and $q=-e^i$, then for any $(t,y,z)\in[0,T]\times\mathbb{R}^n\times\mathbb{R}^{n\times d}$, we have $g^1_i(t,y,z)=g^2_i(t,y,z)$, $P$-a.s., where $g_i^1$ and $g_i^2$ mean the $i$th component of $g^1$ and $g^2$ respectively.
\end{corol}
Consider some special cases of total order $\succsim^q$, using Theorem \ref{th21} or Theorem \ref{th22}, we have the following remarks hold.
\begin{remark}\label{re1} Let $n=1,q=1$. For this case, $x\succsim^1 y$ means $x\geq y$. Then, (i) is equivalent to $g^1(t,y,z)\geq g^2(t,y,z),$ $P$-a.s..
\end{remark}
This coincides with $1$-dimensional comparison theorem established in \cite{f}.\\
{\it Proof.}\,\, ``$\Rightarrow$'': Immediately from Theorem \ref{th22}.\\
``$\Leftarrow$'': We only need to prove that (ii) holds. When $y\geq 0$, (ii) obviously holds. So we only consider the case $y<0$. The left side of inequality (\ref{eq22}) equals
\begin{eqnarray*}
4y[g^1(t,y',z)-g^2(t,y',z')]&\leq&4y[g^2(t,y',z)-g^2(t,y',z')]\\
&\leq&{2\over{\mu^2}}|g^2(t,y',z)-g^2(t,y',z')|^2+2\mu^2y^2\\
&\leq&2|z-z'|^2+2\mu^2y^2,\quad P\hbox{-a.s.}.
\end{eqnarray*}
Letting $C=2\mu^2$ implies that (ii) holds. The proof of Remark \ref{re1} is completed. $\hfill{} \Box$
\begin{remark}\label{re2} Let $q=e^i$. For this case, $y^1\succsim^qy^2$ means $y_i^1\geq y_i^2$, where $y_i^1$ and $y_i^2$ represent the $i$th component of $y^1$ and $y^2$ respectively. Then the following statements are equivalent:
\\$\bullet\ $ (iv) For any $u\in[0,T]$, $\xi^1,\xi^2\in L^2(\Omega,\mathcal{F}_{u},P;\mathbb{R}^n)$ such that $\xi^1_i\geq\xi^2_i $, then the unique solutions $(Y^1,Z^1)$ and $(Y^2,Z^2)$ in $\mathcal{S}^2_{u}\times\mathcal{H}^2_{u}$ to BSDEs (\ref{eq21}) over time interval $[0,u]$ satisfy $(Y_t^1)_i\geq (Y_t^2)_i$, $P$-a.s., for $t\in[0,u];$
\\$\bullet\ $ (v) For any $ t\in[0,T]$, $ (y,z),(y',z')\in\mathbb{R}^n\times\mathbb{R}^{n\times d},$ we have
\begin{equation*}
-4y_i^-\left[g^1_i(t,y+y_i^-e^i+y',z)-g^2_i(t,y',z')\right]\leq2I_{y_i<0}|z_i-z'_i|^2+C(y_i^-)^2,\quad P\hbox{-a.s.},
\end{equation*}
where $C >0$ is a constant independent of $(t, y, z)$.
\end{remark}
The proof is straightforward by taking $q=e^i$ in Theorem \ref{th21}. From Remark \ref{re2}, consider the comparison theorem for the $i$th component in multidimensional BSDEs, we can only set $q=e^i$ in Theorem \ref{th21}. It does not need to consider other components. However, if we use the comparison theorem of multidimensional BSDEs in \cite{hu} to consider the comparison theorem of the $i$th component, we have to do some restrictions on other components.
\begin{example}\label{ex1}
Let $n=2$. For any $t\in[0,T]$, $y=(y_1,y_2)^{\mathbf{T}}\in\mathbb{R}^2$, $z\in\mathbb{R}^{2\times d}$, we have $g^1(t,y,z)=(y_1+y_2, t)^{\mathbf{T}}$, and $g^2(t,y,z)=(y_1+y_2-1, t)^{\mathbf{T}}$. Obviously, $g^1(t,y,z)\succsim ^{e_1}g^2(t,y,z)$.
However, statement (iv) in Remark \ref{re2} does not hold for $i=1$. To prove it, we just suppose that (v) holds for $i=1$. Then for any $ y'=(y'_1,y'_2)^{\mathbf{T}}$, $y=(y_1,-3)^{\mathbf{T}}$, $z=z'$,
where $y_1$ being any negative number. We get $C>-{8\over{y_1}}$, which contradicts the condition that $C$ is a positive constant independent of $y$. Thus, (v) dose not hold for $i=1$, so does (iv).
\end{example}
\begin{remark}\label{re3}
Let $g^1=g^2=g$, and $q=e^i$. Then, (iv) is equivalent to that for any $t\in[0,T]$, the $i$th component of $g$ denoted by $g_i$ depends only on $y_i$ and $z_i$, $P$-a.s..
\end{remark}
{\it Proof.}\,\, We only need to show (v) is equivalent that for any $t$, $g_i$ depends only on $y_i,z_i.$
\\``$\Rightarrow$'': In (v), choose $y=-{1\over k}e^i$, where $k$ is an arbitrary number in $\mathbb{N}^*$. Then for any $y'\in \mathbb{R}^n$, and any $z,z'\in\mathbb{R}^{n\times d}$ such that $z_i=z'_i$, we have $$-4{1\over k}[g_i(t,y',z)]-g_i(t,y',z')\leq C{1\over {k^2}},\quad P\hbox{-a.s.}.$$
Let $k\rightarrow\infty$, we deduce that $g_i(t,y',z)\geq g_i(t,y',z')$, $P$-a.s.. Similarly, we obtain $g_i(t,y',z')\geq g_i(t,y',z)$, $P$-a.s.. Hence, $g_i(t,y',z')=g_i(t,y',z)$, $P$-a.s.. Therefore, for any $(t,y)$, $g_i$ depends only on $z_i$, $P$-a.s..
\par
For any $\bar{y}$ such that $\bar{y}_i=0$, let $y=\bar{y}-\epsilon e^i,\epsilon>0.$ Then for any $y'$, $z=z'$ in (v), we get
$$-4\epsilon [g_i(t,\bar{y}+y',z)-g_i(t,y',z)]\leq C\epsilon^2,\quad P\hbox{-a.s.}.$$
Letting $\epsilon \rightarrow 0$, we deduce that $g_i(t,\bar{y}+y',z)\geq g_i(t,y',z)$, $P$-a.s.. Noticing the property of $\bar{y}$, we also get $g_i(t,\bar{y}+y',z)\leq g_i(t,y',z)$, $P$-a.s.. Therefore, $g_i(t,\bar{y}+y',z)=g_i(t,y',z),$ $P$-a.s., that is, for any $(t,z)$, $g_i$ depends only on $y_i$, $P$-a.s..
\par
From the arguments above it follows that for any $t$, $g_i$ depends only on $y_i$ and $z_i$, $P$-a.s..
\\``$\Leftarrow$'': It's clearly that (v) is always true for $y_i\geq 0$. Now we consider $y_i<0$, and have
\begin{eqnarray*}
4y_i[g_i(t,y-y_ie^i+y',z)-g_i(t,y',z')]&=&4y_i[g_i(t,y'_i,z_i)-g_i(t,y'_i,z'_i)]\\
&\leq&{2\over{\mu^2}}|g_i(t,y'_i,z_i)-g_i(t,y'_i,z'_i)|^2+2\mu^2y_i^2\\
&\leq&2|z_i-z_i'|^2+2\mu^2y_i^2,\quad P\hbox{-a.s.}.
\end{eqnarray*}
Thus, letting $C=2\mu^2$, (v) immediately holds. The proof of Remark \ref{re3} is completed. $\hfill{} \Box$
\begin{remark}\label{re4} If (iv) holds for all $i\in\{1,\cdots,n\}$, we have\\
$\bullet\ $ (vi) for any $u\in[0,T]$, $\xi^1,\xi^2\in L^2(\Omega,\mathcal{F}_{u},P;\mathbb{R}^n)$ such that $\xi^1\geq \xi^2 $, then the unique solutions $(Y^1,Z^1)$ and $(Y^2,Z^2)$ in $\mathcal{S}^2_{u}\times\mathcal{H}^2_{u}$ to BSDEs (\ref{eq21}) over time interval $[0,u]$ satisfy
$$Y_t^1\geq Y_t^2, \quad P\hbox{-a.s.},\quad \forall t\in[0,u].$$
\end{remark}
The proof of Remark \ref{re4} is obviously. However, the converse of Remark \ref{re4} is not true, that is (iv) can not be deduced from (vi). We use the following example to illustrate this.
\begin{example}\label{ex2}
Let $n=2$. Then for any $ t\in[0,T]$, $z\in\mathbb{R}^{2\times d}$, $y=(y_1,y_2)^{\mathbf{T}}\in\mathbb{R}^2$, $g^1(t,y,z)=g^2(t,y,z)=(y_1+y_2,|z_2|)^{\mathbf{T}}$, due to Theorem 2.2 in \cite{hu}, (vi) holds.
However, from Remark \ref{re3}, (iv) does not hold for $i=1$.
\end{example}
\begin{remark}\label{re5} Let $q=({1\over{\sqrt{n}}},\cdots,{1\over{\sqrt{n}}})^{\mathbf{T}}={1\over{\sqrt{n}}}\mathbf{1}$, where $\mathbf{1}=(1,\cdots,1)^{\mathbf{T}}\in\mathbb{R}^n$. For this case, it is seen that the statement (i) holds for $q$ is equivalent to the following:
\begin{eqnarray*}
&&-4(\sum_{j=1}^ny_j)^-\left[\sum_{i=1}^n\left(g_i^1(y+{\mathbf{1}\over n}(\sum_{j=1}^ny_j)^-+y',z)-g_i^2(y',z')\right)\right]\\
&&\leq 2I_{\sum_{j=1}^ny_j<0}\sum_{i,j=1}^n(z_i-z'_i)^{\mathbf{T}}(z_i-z'_j)+C[(\sum_{j=1}^ny_j)^-]^2,\quad P\hbox{-a.s.}.
\end{eqnarray*}
\end{remark}
{\bf Acknowledgement.}\hspace{0.4cm} The authors would like to thank the editors and the referees for their careful reading and useful suggestions. And also, the authors are very grateful to Dr. Yu Pei, the Professor of University of Western Ontario, who read the paper after review and corrected uncountable misuses of the English language.

\end{document}